\documentclass{amsart}

\usepackage{amssymb}
\usepackage{bbm}
\usepackage{float}

\makeatletter
\def\hsmash{\relax % \relax, in case this comes first in \halign
  \ifmmode\def\next{\mathpalette\mathhsm@sh}\else\let\next\makehsm@sh
  \fi\next}
\def\makehsm@sh#1{\setbox\z@\hbox{#1}\finhsm@sh}
\def\mathhsm@sh#1#2{\setbox\z@\hbox{$\m@th#1{#2}$}\finhsm@sh}
\def\finhsm@sh{\wd\z@\z@ \box\z@}
\makeatother

\newcommand{\br}{ }
\newcommand{\brr}{, }

\newtheorem{fac}{Fact}[section]
\newtheorem{lem}[fac]{Lemma}
\newtheorem{prop}[fac]{Proposition}
\newtheorem{theo}[fac]{Theorem}
\newtheorem{coro}[fac]{Corollary}

\theoremstyle{definition}
\newtheorem{ttt}[fac]{}
\newtheorem{algo}[fac]{Algorithm}
\newtheorem{obs}[fac]{Observations}

\theoremstyle{remark}
\newtheorem{rem}[fac]{Remark}
\newtheorem{rems}[fac]{Remarks}
\newtheorem{ex}[fac]{Example}
\newtheorem{nota}[fac]{Notation}

\def\pmod#1{\nobreak\ifinner\mkern8mu\else\mkern8mu\fi (\text{\rmfamily\upshape mod}\,\,#1)}

\newcommand{\bP}{\mathop{\text{\bf P}}\nolimits}
\newcommand{\bA}{\mathop{\text{\bf A}}\nolimits}

\newcommand{\Br}{\mathop{\text{\rm Br}}\nolimits}
\newcommand{\Kum}{\mathop{\text{\rm Kum}}\nolimits}
\newcommand{\Pic}{\mathop{\text{\rm Pic}}\nolimits}
\newcommand{\NS}{\mathop{\text{\rm NS}}\nolimits}
\newcommand{\Gal}{\mathop{\text{\rm Gal}}\nolimits}

\newcommand{\inv}{\mathop{\text{\rm inv}}\nolimits}
\newcommand{\sep}{\mathop{\text{\rm sep}}\nolimits}
\newcommand{\Hom}{\mathop{\text{\rm Hom}}\nolimits}
\newcommand{\im}{\mathop{\text{\rm im}}\nolimits}
\newcommand{\cha}{\mathop{\text{\rm char}}\nolimits}

\newcommand{\ev}{\mathop{\text{\rm ev}}\nolimits}

\newcommand{\bbA}{{\mathbbm A}}
\newcommand{\bbF}{{\mathbbm F}}
\newcommand{\bbQ}{{\mathbbm Q}}
\newcommand{\bbR}{{\mathbbm R}}
\newcommand{\bbZ}{{\mathbbm Z}}

\newcommand{\good}{{{\text{\rm good}}}}
\newcommand{\et}{{{\text{\rm {\'e}t}}}}

\renewcommand{\atop}[2]{\genfrac{}{}{0pt}{}{#1}{#2}}

\newcounter{abc}
\newenvironment{abc}{\begin{list}{\rm \alph{abc}) }%
{\usecounter{abc} \leftmargin=0.0pt \labelsep=0.0pt %
\listparindent=0.0pt \labelwidth=0.0pt \parsep=\smallskipamount %
\itemsep=0.0pt \topsep=0.0pt \partopsep=\smallskipamount}}{\end{list}}

\newcounter{iii}
\newenvironment{iii}{\begin{list}{\rm \roman{iii}) }%
{\usecounter{iii} \leftmargin=0.0pt \labelsep=0.0pt %
\listparindent=0.0pt \labelwidth=0.0pt \parsep=\smallskipamount%
 \itemsep=0.0pt \topsep=0.0pt \partopsep=\smallskipamount}}{\end{list}}

\newcounter{III}

\def\rightend#1#2{{%
 \leavevmode\nobreak\hskip .5em plus 1fil
 \penalty600 \hskip 0pt plus -1filll
 \vadjust{}\nobreak\hskip 0pt plus 1filll%
 #1\parfillskip=#2\relax \par}}

\def\eop{\ifmmode\rule[-22pt]{0pt}{1pt}\ifinner\tag*{$\square$}\else\eqno{\square}\fi\else\rightend{$\square$}{0pt}\fi}

\title[Transcendental Brauer-Manin obstruction]{Experiments with the transcendental Brauer-Manin obstruction}

\begin{document}

\author{Andreas-Stephan Elsenhans}

\address{Mathematisches Institut\\ Universit\"at Bayreuth\\ Universit\"ats\-stra\ss e 30\\ D-95440 Bay\-reuth\\ Germany}
\email{Stephan.Elsenhans@uni-bayreuth.de}
\urladdr{http://www.staff\!.\!uni-bayreuth.de/$\sim$btm216}

\author{J\"org Jahnel}

\address{D\'epartement Mathematik\\ Universit\"at Siegen\\ Walter-Flex-Stra\ss e~3\\ D-57068 Sie\-gen\\ Germany}
\email{jahnel@mathematik.uni-siegen.de}
\urladdr{http://www.uni-math.gwdg.de/jahnel}

%\thanks{}

\date{May~13,~2012.}

\keywords{Kummer surface, Weak approximation, Transcendental Brauer-Manin obstruction, Multivariate paging}

\subjclass[2010]{Primary 11D41; Secondary 11Y50, 11G35, 14J28}

\begin{abstract}
We~report on our experiments and theoretical investigations concerning weak approximation and the transcendental Brauer-Manin obstruction for special Kummer surfaces.
\end{abstract}

\maketitle

\section{Introduction}

\begin{ttt}
Consider a geometrically integral, projective variety
$S$
over the field
$\bbQ$
of rational~numbers. It~is said that
$S$
fulfills weak approximation when the following is~true. For~every finite
set~$\{p_1,\ldots,p_l\}$
of prime numbers and every vector
$$(x_0, x_1, \ldots, x_l) \in S(\bbR) \times S(\bbQ_{p_1}) \times \cdots \times S(\bbQ_{p_l}) \, ,$$
there exists a sequence of
\mbox{$\bbQ$-rational}
points that simultaneously converges to
$x_i$
in the
\mbox{$p_i$-adic}
topology for
$i = 1, \ldots, l$
and to
$x_0$
with respect to the real~topology. In~a more formal language, this means that the set
$S(\bbQ)$
of the rational points
on~$S$
is dense in the set
$S(\bbA_\bbQ)$
of all adelic~points.

Even~for Fano varieties, which are generally expected to have many rational points, weak approximation is not always fulfilled. Well-known counterexamples are due to Sir Peter Swinnerton-Dyer~\cite{SD}, L.\,J.~Mordell~\cite{Mo}, J.\,W.\,S.~Cassels and M.\,J.\,T.~Guy~\cite{CG}, and many~others.

For varieties of intermediate type,
e.g.~$K3$~surfaces,
the situation is yet more~obscure. In~fact, to prove the much weaker statement that
$\#S(\bbQ) = \infty$
is usually a formidable task in its own~\cite{LKL,Ka}. It~seems therefore that proving weak approximation, even for a single
$K3$~surface,
is presently out of reach and that experiments are asked~for.
\end{ttt}

\begin{ttt}
To test weak approximation experimentally is, however, an ill-posed problem, at least from the strictly formal point of~view. The~reason is that weak approximation is not a finite~phenomenon. It~is strongly infinite in~nature.

An~interesting situation occurs when a certain ``obstruction'' is responsible for the failure of weak~approximation. This~means that
$S(\bbQ_p)$
breaks somehow regularly into open-closed subsets, each of which behaves uniformly as far as approximation by
$\bbQ$-rational
points is~concerned.
As~$S(\bbQ_p)$
is compact, it is clear that finitely many subsets
$U_1, \ldots, U_k \subset S(\bbQ_p)$
will~suffice. When such a behaviour appears, we speak of a {\em colouring\/} and call the subsets the {\em colours\/}
of~$S(\bbQ_p)$.
\end{ttt}

\begin{ttt}
\label{reg}
It is well-known that a class
$\alpha \in \Br(S)$
in the Grothendieck-Brauer group
of~$S$
induces such a~colouring.
For~a
point~$x \in S(\bbQ_p)$,
its colour is obtained as
$\inv_{\bbQ_p} (\alpha|x) \in \bbQ/\bbZ$.
If~$\alpha$
is of
order~$N$
then not more than
$N$~colours
may~occur.

As~a result, a failure of weak approximation may~appear. Indeed,~for a
\mbox{$\bbQ$-rational}
point, one must have
$\sum_p \inv_{\bbQ_p} (\alpha|x) = 0$,
but the same need not be true for an adelic~point. This~phenomenon is called the Brauer-Manin obstruction~\cite{Ma}.

There~is a canonical filtration
on~$\Br(S)$,
which causes a distinction between {\em algebraic\/} and {\em transcendental\/} Brauer~classes. Correspondingly,~there are the algebraic and the transcendental Brauer-Manin~obstructions.

The~algebraic Brauer-Manin obstruction is rather well understood. At~least on
$S(\bbQ_p)_\good \subseteq S(\bbQ_p)$,
the
\mbox{$p$-adic}
points with good reduction, it yields extremely regular colourings~\cite{CKS,CS,EJ3}. For~example, a colouring by two colours is possible only when there is an unramified two-sheeted covering
$\pi\colon X \to S(\bbQ_p)_\good$.
The~two colours are then given by the subsets
$\{ x \in S(\bbQ_p) \mid \pi^{-1}(x) = \emptyset \}$
and
$\{ x \in S(\bbQ_p) \mid \#\pi^{-1}(x) = 2 \}$.

Explicit computations of the algebraic Brauer-Manin obstruction have been done for many classes of~varieties. Most~of the examples were~Fano. For~instance, we gave a systematic treatment of the (algebraic) Brauer-Manin obstruction for cubic surfaces in~\cite{EJ2,EJ4}.
Concerning~$K3$
surfaces, computations for diagonal quartic surfaces are due to M.\,Bright~\cite{Br}. Further,~it is known that there is no algebraic Brauer-Manin obstruction on a generic Kummer surface as well as on the generic case of a Kummer surface associated to the product of two elliptic curves~\cite[Proposition~1.4.ii)]{SZ}.
\end{ttt}

\begin{ttt}
The~transcendental Brauer-Manin obstruction is much less understood and seems to be by far more difficult, at least at~present. The~historically first example of a variety, where weak approximation is violated due to a transcendental Brauer class, was constructed by D.\,Harari~\cite{Ha}.

Concerning~$K3$~surfaces,
the available literature is still rather~small. The~interested reader is encouraged to consult the articles~\cite{Wi}, \cite{SSD}, \cite{Ie}, \cite{ISZ}, \cite{Pr}, \cite{HVV}, and~\cite{HV}, at least in order to recognize the enormous efforts made by the~authors. For~example, the whole Ph.D.\,thesis of Th.\,Preu is devoted to the computation of the transcendental Brauer-Manin obstruction for single diagonal quartic~surface.

An~exceptional case, which seems to be a bit more accessible, is provided by the Kummer surfaces
$S := \Kum(E \times E')$
for two elliptic curves
$E$
and~$E'$.
Here,~the Brauer group, which is typically purely transcendental, was described in detail by A.\,N.\ Skoro\-bogatov and Yu.\,G.\ Zarhin in~\cite{SZ}.
\end{ttt}

\begin{ttt}
For this reason, in the present article, we will deal with Kummer surfaces, defined
over~$\bbQ$,
of this particular~type. To~keep the theory simple, we will restrict ourselves to the case that both curves have their full
\mbox{$2$-torsion}
defined over the base~field. We~may start with equations of the form
$E\colon y^2 = x(x-a)(x-b)$
and
$E'\colon y^2 = x(x-a')(x-b')$,
for~$a,b,a',b' \in \bbQ$.
Then~$S := \Kum(E \times E')$
is a double cover of
$\smash{\bP^1 \times \bP^1}$,
an affine chart of which is given by the~equation
\begin{equation}
\label{sext}
z^2 = x(x-a)(x-b) u(u-a')(u-b') \, .
\end{equation}
The~goal of the article is to report on our experiments and theoretical investigations concerning weak approximation and the transcendental Brauer-Manin obstruction for Kummer surfaces of this particular~type.
\end{ttt}

\begin{rem}
To~be precise, equation~(\ref{sext}) defines a model of the Kummer surface with 16 singular points of
type~$A_1$.
In~the minimal regular model, the singularities are replaced by projective~lines.
As~$\smash{\Br(\bP^1_k) = \Br(k)}$,
the evaluation of a Brauer class on a projective line is automatically~constant. Thus,~we may work as well with the singular~model.
\end{rem}

\begin{ttt}
{\bf The results.}
Among the Kummer surfaces of type~(\ref{sext}) for
$|a|,|b|,|a'|$,
$|b'| \leq 200$,
we determined all those, for which there is a transcendental Brauer-Manin obstruction arising from a
\mbox{$2$-torsion}
Brauer~class.

We~found out that there were exactly 3418 surfaces having a nontrivial
\mbox{$2$-torsion}
Brauer~class. In~three cases, this class was~algebraic. Moreover,~we identified the adelic subsets of the surfaces where the Brauer class gives no~obstruction. On~only six of the surfaces, it happened that no adelic point was~excluded. 

On~the other hand, we developed a memory-friendly point searching algorithm for Kummer surfaces of the form~above. The~sets of
$\bbQ$-rational
points found turned out to be compatible with the idea that the Brauer-Manin obstruction might be the only obstruction to weak~approximation.
\end{ttt}

\section{The transcendental Brauer group}

\subsubsection*{Generalities.}

\begin{nota}
For
$A$
an abelian group, we will write
$A_2$
to denote its
$2$-torsion~part.
\end{nota}

\begin{ttt}
The cohomological Grothendieck-Brauer group of an algebraic variety
$S$
over a
field~$k$
is equipped with a canonical three-step~filtration, defined by the Hochschild-Serre spectral~sequence.

\begin{iii}
\item
$\Br_0(S) \subseteq \Br(S)$
is the image
of~$\Br(k)$
under the natural~map. At~least when
$S$
has a
\mbox{$k$-rational}~point,
$\Br_0(S) \cong \Br(k)$.
For~$k$
a number field, the existence of an adelic point~suffices.
$\Br_0(S)$~does
not contribute to the Brauer-Manin~obstruction.
\item
The~quotient
$\Br_1(S)/\Br_0(S)$
is isomorphic to
$H^1(\Gal(k^{\sep}/k), \Pic(S_{k^{\sep}}))$.
This~subquotient is called the algebraic part of the Brauer~group.
For~$k$
a number field, it is responsible for the algebraic Brauer-Manin~obstruction.
\item
Finally,~$\Br(S)/\Br_1(S)$
injects into
$\Br(S_{k^{\sep}})$.
This~quotient is called the transcendental part of the Brauer~group. Nevertheless,~every Brauer class that is not algebraic is usually said to be~transcendental.
For~$k$
a number field, the corresponding obstruction is a transcendental Brauer-Manin~obstruction.
\end{iii}
\end{ttt}

\begin{ttt}
\label{Skoro}
For~$S$,
the Kummer surface corresponding to the product of two elliptic curves, the Brauer group
of~$S$
is well understood due to the work~\cite{SZ} of A.\,N.~Skorobogatov and Yu.\,G.~Zarhin. For~us, the following result will be~sufficient.\medskip

\noindent
{\bf Proposition}
(Skorobogatov/Zarhin)%
{\bf .}
{\em
Let\/
$E \colon y^2 = x(x-a)(x-b)$
and\/
$E' \colon v^2 = u(u-a')(u-b')$
be two elliptic curves over a
field\/~$k$
of characteristic~zero. Suppose~that their\/
\mbox{$2$-torsion}
points are defined
over\/~$k$
and that\/
$\smash{E_{\overline{k}}}$
and\/
$\smash{E'_{\overline{k}}}$
are not isogenous to each~other.\smallskip

\noindent
Further,~let\/
$S := \Kum(E \times E')$
be the corresponding Kummer~surface. Then
$$\Br(S)_2/\Br(k)_2 = \im(\Br(S)_2 \to \Br(S_{\overline{k}})_2) \cong \ker(\mu\colon \bbF_{\!2}^4 \to (k^*/k^{*2})^4) \, ,$$
where\/
$\mu$
is given by the matrix
\begin{equation}
\label{szmat}
M_{aba'b'} :=
\left(
\begin{array}{cccc}
  1  &       ab  &   a'b' &      -aa' \\
 ab  &        1  &    aa' & a'(a'-b') \\
a'b' &       aa' &      1 &    a(a-b) \\
-aa' & a'(a'-b') & a(a-b) &         1
\end{array}
\right) .
\end{equation}
}%
{\bf Proof.}
The equality on the left hand side expresses the absence of algebraic Brauer classes, which is shown in \cite[Proposition~3.5.i]{SZ}. The~isomorphism on the right is established in \cite[Proposition~3.5.ii and~iii, together with Lemma~3.6]{SZ}. The~reader might want to compare~\cite[Proposition~3.7]{SZ}.
\eop
\end{ttt}

\begin{ttt}
Consider the case where
$k$
is algebraically~closed. Then,~induced by the Kummer sequence, there is the short exact~sequence
$$0 \rightarrow \Pic(S)/2\Pic(S) \rightarrow H^2_\et(S, \mu_2) \rightarrow \Br(S)_2 \rightarrow 0 \, .$$
We~have
$\dim_{\bbF_{\!2}} \Pic(S)/2\Pic(S) = 16 + \dim_{\bbF_{\!2}} \NS(E \times E') / 2\NS(E \times E') = 18$
and
$\dim_{\bbF_{\!2}} H^2_\et(S, \mu_2) = 22$.
This explains why
$\Br(S)_2 \cong \bbF_{\!2}^4$.
More~canonically, there are~isomorphisms
$$\Br(S)_2 \cong H^2_\et(E \times E', \mu_2) / (H^2_\et(E, \mu_2) \oplus H^2_\et(E'\!, \mu_2))\cong \Hom(E[2], E'[2]) \, .$$
\end{ttt}

\begin{rem}
For~$k$
any field of characteristic zero, the assumption that the
\mbox{$2$-torsion}
points are defined
over~$k$
therefore implies that
$\Gal(\overline{k}/k)$
operates trivially
on~$\Br(S_{\overline{k}})_2$.
%I.e.,~$\smash{\Br(S_{\overline{k}})_2^{\Gal(\overline{k}/k)} \cong \bbF_{\!2}^4}$.
We~see explicitly that
$\smash{\Br(S)_2 \subsetneqq \Br(S_{\overline{k}})_2^{\Gal(\overline{k}/k)} \cong \bbF_{\!2}^4}$,
in~general.
\end{rem}

\begin{ttt}
Assume~that
$k$
is algebraically~closed.
For~two rational functions
$f,g \in k(S)$,
we denote by
$(f,g)$
the quaternion algebra
$$k(S)\{I,J\}/(I^2 - f, J^2 - g, IJ + JI)$$
over~$k(S)$.
Cohomologically,~$f$
and~$g$
define classes
in~$H^1(\Gal(\overline{k(S)}/k(S)), \mu_2)$
via the Kummer~sequence. The~Brauer class
of~$(f,g)$
is the cup product in
\begin{eqnarray*}
H^2(\Gal(\overline{k(S)}/k(S)), \mu_2^{\otimes 2}) & = & H^2(\Gal(\overline{k(S)}/k(S)), \mu_2) \\
 & \subseteq & H^2(\Gal(\overline{k(S)}/k(S)), \overline{k(S)}^*)
\end{eqnarray*}
of these two~classes. The~symbol
$(.\,,.)$
is thus bilinear and~symmetric.
\end{ttt}

\begin{fac}
\label{expl}
Let\/~$k$
be an algebraically closed field,
$\cha k = 0$,
$a,b,a',b' \in k$,
and\/
$S$
be as in Proposition~\ref{Skoro}. Then,~in terms of the canonical injection\/
$\Br(S) \hookrightarrow \Br(k(S))$,
a basis
of\/~$\Br(S)_2$
is given by the four quaternion~algebras
$$A_{\mu,\nu} := ((x-\mu)(x-b), (u-\nu)(u-b')), \qquad \mu = 0,a, \;\nu = 0,a'.$$
Here,~the standard vectors
in\/~$\bbF_{\!2}^4$
correspond to these four~algebras. More~precisely,
$e_1$
corresponds to\/
$A_{a,a'}$,
$e_2$
to\/
$A_{a,0}$,
$e_3$
to\/
$A_{0,a'}$,
and\/
$e_4$
to\/~$A_{0,0}$.\smallskip

\noindent
{\bf Proof.}
{\em
This~is \cite[Lemma~3.6]{SZ} together with \cite[formula~(20)]{SZ}.
}
\eop
\end{fac}

\begin{rem}
Using bilinearity, for nine of the 15 non-trivial classes, we find a description as a single quaternion algebra similar to the type~above. For~the six classes corresponding to the vectors
$(1,0,0,1)$,
$(0,1,1,0)$,
$(1,1,1,0)$,
$(1,1,0,1)$,
$(1,0,1,1)$,
and
$(0,1,1,1)$,
we need at least two such~algebras.
\end{rem}

\begin{obs}[Isomorphy, Twisting]
\label{isotw}
\begin{iii}
\item
We~may replace
$(a,b)$
by
$(-a, b-a)$
or~$(-b, a-b)$
without
changing~$S$,
and similarly
for~$(a',b')$.
Indeed,~this simply means to apply the translations
$\smash{\bA^1_k \to \bA^1_k}$,
$x \mapsto x - \mu$,
for~$\mu = a, b$.
\item
It~is also possible to replace
$(a,b,a',b')$
by
$(\lambda^2 a, \lambda^2  b,a',b')$
or
$(\lambda a, \lambda b, \lambda a', \lambda b')$
for~$\lambda \in k$.
The~reason is that the twist
$E^{(\lambda)}\colon \lambda y^2 = x(x-a)(x-b)$
is isomorphic to the elliptic curve, given by
$Y^2 = X(X - \lambda a)(X - \lambda b)$.
\end{iii}
\end{obs}

\begin{ttt}[The isogenous case]
When
$E_{\overline{k}}$
and~$E'_{\overline{k}}$
are isogenous, only minor modifications~occur. The~isogeny causes
$\smash{\NS(E_{\overline{k}} \times E'_{\overline{k}}) / 2\NS(E_{\overline{k}} \times E'_{\overline{k}})}$
to have dimension higher than~two. Hence,~the homomorphism
$\bbF_{\!2}^4 \cong \Hom(E[2], E'[2]) \to \Br(S_{\overline{k}})_2$
is only a surjection, not a~bijection.

Over~a non-algebraically closed field, the situation is as~follows.
If~$E$
and~$E'$
are isogenous
over~$k$
then
$\dim_{\bbF_{\!2}} \!\Pic(S)/2\Pic(S) > 16 + 2 = 18$.
As~the additional generator evaluates trivially, it will be found
in~$\ker M_{aba'b'}$~\cite[Lemma~3.6]{SZ}.
Thus,~the homomorphism
$\ker M_{aba'b'} \twoheadrightarrow \Br(S)_2/\Br(k)_2$
has a non-trivial~kernel.

An~isogeny defined over a proper field extension
$l/k$
causes the same effect
over~$l$,
but not
over~$k$.
As~$\Pic(S)/2\Pic(S) \subsetneq \Pic(S_l)/2\Pic(S_l)$,
it may, however, happen that a Brauer class is annihilated by the
extension~$l/k$.
I.e.,~that a vector in
$\ker M_{aba'b'}$
describes an {\em algebraic\/} Brauer~class. By~the Hochschild-Serre spectral sequence, we have
$H^2_\et(S,\mu_2)/\Br(k)_2 \subseteq H^2_\et(S_{\overline{k}},\mu_2)$.
Hence,~there are no other algebraic
\mbox{$2$-torsion}
Brauer classes than~these.
\end{ttt}

\subsubsection*{The transcendental Brauer-Manin obstruction}

\begin{lem}
Let\/~$k$
be a local field of characteristic~zero. For~two elliptic curves\/
$E \colon y^2 = x(x-a)(x-b)$
and\/
$E' \colon v^2 = u(u-a')(u-b')$
over\/~$k$
with\/
\mbox{$k$-rational\/}
\mbox{$2$-torsion},
consider\/
$S := \Kum(E \times E')$,
given explicitly~by
$$z^2 = x(x-a)(x-b) u(u-a')(u-b') \, .$$
Let\/~$\alpha \in \Br(S)$
be a Brauer class, represented by an Azumaya algebra
over\/~$k(S)$
of the
type\/~$\bigotimes_i A_{\mu_i,\nu_i}$.\smallskip

\noindent
Then~the local evaluation map\/
$\ev_\alpha\colon S(k) \to \frac12\bbZ/\bbZ$
is given by
$$(x,u;z) \mapsto \ev_\alpha((x,u;z)) = \sum_i ((x-\mu_i)(x-b), (u-\nu_i)(u-b'))_k \, .$$
Here,~$(.\,,.)_k$
denotes the
$k$-Hilbert
symbol\/~\cite[Ch.\,1, \S6]{BS}%
%followed by the isomorphism
%$(\{\pm1\}, \;\cdot\;) \to (\frac12\bbZ/\bbZ, +)$
.\medskip

\noindent
{\bf Proof.}
{\em
By~definition,
$\ev_\alpha((x,u;z)) = \inv(\alpha|_{(x,u;z)})$.
Further,~$\alpha|_{(x,u;z)}$
is the Azumaya algebra
$\bigotimes_i ((x-\mu_i)(x-b), (u-\nu_i)(u-b'))$
over~$k$.
Now~observe that the quaternion algebra
$(s,t)$
splits if and only if
$t$
is a norm
from~$k(\sqrt{s})$.
This~is tested by the norm residue symbol
$(t, k(\sqrt{s})/k)$,
which agrees with the classical Hilbert
symbol~$(s,t)_k$.
}
\eop
\end{lem}

\begin{rems}
\begin{iii}
\item
For~us, the Hilbert symbol takes values
in~$(\frac12\bbZ/\bbZ, +)$.
This~differs from the classical setting, where the values are taken
in~$(\{\pm1\}, \,\cdot\,)$.
\item
According~to Proposition~\ref{Skoro},
$\Br(S)_2/\Br(k)_2 \subseteq \bbF_{\!2}^4$.
Further,~by Fact~\ref{expl}, we have an explicit basis, which is given by Azumaya~algebras. I.e.,~for each class
in~$\Br(S)_2/\Br(k)_2$,
we chose a lift
to~$\Br(S)_2$.

For~$k$
a local field, this lift is normalized such that 
$\ev_\alpha((\infty,\infty,\;.\;)) = 0$.
Indeed,~for
$x$
close
to~$\infty$
in~$k$,
$(x-\mu)(x-b)$
is automatically a~square.
\end{iii}
\end{rems}

\begin{ttt}
\label{appr}
The~evaluation map is constant near the singular~points.\smallskip

\noindent
{\bf Lemma.}
{\em
Let\/~$p > 2$
be a prime number and\/
$a,b,a',b' \in \bbZ_p$
be such that\/
$E\colon y^2 = x(x-a)(x-b)$
and\/
$E'\colon v^2 = u(u-a')(u-b')$
are elliptic curves, not isogenous to each~other. Put
$$l := \max (\nu_p(a), \nu_p(b), \nu_p(a-b), \nu_p(a'), \nu_p(b'), \nu_p(a'-b')) \, .$$
Consider~the
surface\/~$S$
over\/~$\bbQ_p$,
given by\/
$z^2 = x(x-a)(x-b) u(u-a')(u-b')$.
Then,~for every
$\alpha \in \Br(S)_2$,
the evaluation map\/
$S(\bbQ_p) \to \bbQ/\bbZ$
is constant on the~subset
\begin{align*}
T := \{(x,u;z) \in S(\bbQ_p) \mid \nu_p(x) < 0 \text{~or~} \nu_p(u) < 0 \text{~or~}& \\[-1mm]
x \equiv \mu, u \equiv \nu \pmod {p^{l+1}}&, \,\mu = 0,a,b, \,\nu = 0,a',b' \} \, .
\end{align*}%
}%
{\bf Proof.}
Consider~the Hilbert symbol
$((x-a)(x-b), (u-a')(u-b'))_p$,~first.
Using~the equation of the surface, we see~that
\begin{equation}
\label{hilb}
((x-a)(x-b), (u-a')(u-b'))_p = ((x-a)(x-b), -xu)_p = (-xu, (u-a')(u-b'))_p .
\end{equation}
Let~us distinguish three~cases. In~all cases, we observe that a Hilbert symbol is clearly zero, when at least one of its entries is a~square.\smallskip

\noindent
{\em First case.}
Negative~valuation.

\noindent
It~is clear that
$\nu_p(x) < 0$
implies
$(x-a)(x-b)$
is a square and that
$\nu_p(u) < 0$
implies that
$(u-a')(u-b')$
is a~square.\smallskip

\noindent
{\em Second case.}
$x \equiv \mu, u \equiv \nu \pmod {p^{l+1}}$
for
$(\mu,\nu) = (0,\nu)$
or~$(\mu,0)$.

\noindent\looseness-1
First,~let
$x \equiv 0 \pmod {p^{l+1}}$.
Then~$(x-a)(x-b) \equiv ab \pmod {p^{l+1}}$.
As
$\nu_p(ab) = \nu_p(a) + \nu_p(b) = \max(\nu_p(a), \nu_p(b)) \leq l$,
both numbers belong to the same square~class.\smallskip

\noindent
Analogously,~$u \equiv 0 \pmod {p^{l+1}}$
implies
$(u-a')(u-b') \equiv a'b' \pmod {p^{l+1}}$
such that
$(u-a')(u-b')$
is in the square class
of~$a'b'$.\smallskip

\noindent
{\em Third case.}
$x \equiv \mu, u \equiv \nu \pmod {p^{l+1}}$
for
$(\mu,\nu) = (a,a')$,
$(a,b')$,
$(b,a')$,
or~$(b,b')$.

\noindent
Suppose,~for example that
$x \equiv a \pmod {p^{l+1}}$
and
$u \equiv a' \pmod {p^{l+1}}$.
Then, in particular,
$x \equiv a \pmod {p^{\nu(a)+1}}$
and
$u \equiv a' \pmod {p^{\nu(a')+1}}$.
This implies
$(-xu) \equiv (-aa') \pmod {p^{\nu(a)+\nu(a')+1}}$
such that
$(-xu)$
is in the square class
of
$(-aa')$.
The~other cases yield the square classes of
$(-ab')$,
$(-ba')$,
and~$(-bb')$.\smallskip

Consequently,~the evaluation map is constant on the set described if and only if the~vector
$$(1, ab, a'b', -aa')^t \in (\bbQ_p^*/\bbQ_p^{*2})^4$$
is~zero. This~is exactly the first column of the
matrix~$M_{aba'b'}$,
cf.~formula~(\ref{szmat}).

For~the Hilbert symbols
$((x-a)(x-b), u(u-b'))_p$,
$(x(x-b), (u-a')(u-b'))_p$,
and~$(x(x-b), u(u-b'))_p$,
the calculations are completely~analogous. They~lead to the second, third, and fourth column
of~$M_{aba'b'}$.

From~this, we see that, for a combination of Hilbert symbols, the evaluation map is constant on the
set~$T$
if and only if it represents a Brauer~class.
\eop
\end{ttt}

\begin{rem}
For~$p=2$,
the condition has to be strengthened to
$\nu_2(x) < -2$
or
$\nu_2(u) < -2$
or
$x \equiv \mu, u \equiv \nu \pmod {2^{l+3}}$.
The~proof is essentially the~same.
\end{rem}

\begin{prop}
\label{CoSk}
Let\/
$E \colon y^2 = x(x-a)(x-b)$
and\/
$E' \colon v^2 = u(u-a')(u-b')$
be two elliptic curves over a local
field\/~$k$,
not isogenous to each~other. Suppose~that\/
$a,b,a',b' \in k$.
Further,~let\/
$S := \Kum(E \times E')$
be the corresponding Kummer~surface.\smallskip

\noindent
Suppose that either\/
$k = \bbR$
or\/
$k$
is a\/
\mbox{$p$-adic}
field and both\/
$E$
and\/~$E'$
have good~reduction. Then,~for every\/
$\alpha \in \Br(S)_2$,
the evaluation map\/
$\ev_\alpha\colon S(k) \to \bbQ/\bbZ$
is constant.\smallskip

\noindent
{\bf Proof.}
{\em
{\em First case.}
$k = \bbQ_p$.

\noindent
This~assertion is a particular case of a very general result~\cite[Proposition~2.4]{CS}, due to \mbox{J.-L.}~Colliot-Th\'el\`ene and A.\,N.~Skorobogatov.
%Let~us nevertheless sketch an argument, which is specific for the present situation but more~elementary.
Using~Lemma~\ref{appr} and the elementary properties of the Hilbert symbol, one may as well provide an elementary argument that is specific for the present~situation.
%
%Assume~without loss of generality that
%$a,b,a',b' \in \bbZ_p$.
%Potential good reduction implies
%$p \notd a,b$,
%$a \not\equiv b \pmod p$,
%and the analogous conditions for
%$a'$
%and~$b'$.
%In~particular,~$p > 2$.
%
%Again,~we consider the Hilbert symbol
%$((x-a)(x-b), (u-a')(u-b'))_p$,
%first. In~order to be non-zero, for a Hilbert symbol, it is necessary that at least one of its entries is of odd
%\mbox{$p$-adic}~valuation.
%This~enforces~$\smash{x,u \in \bbZ_p}$.
%Further,~two of the nine congruences
%$x \equiv 0,a,b \pmod p$,
%$\smash{u \equiv 0,a',b' \pmod p}$
%must be~fulfilled. Indeed,~otherwise for at least one of the three representations given in~(\ref{hilb}), we see that the Hilbert symbol is~zero. For~the Hilbert symbols
%$((x-a)(x-b), u(u-b'))_p$,
%$(x(x-b), (u-a')(u-b'))_p$,
%and~$(x(x-b), u(u-b'))_p$,
%the argument is~analogous.
%
%Consequently,~for a Brauer class that is normalized
%to~$\ev_\alpha((\infty,\infty,\;.\;)) = 0$,
%the evaluation map
%$\ev_\alpha\colon S(\bbQ_p) \to \bbQ/\bbZ$
%is zero, except possibly for the points reducing to one of the~singularities. But~for these, Lemma~\ref{appr} is~applicable. Observe~that because of the good reduction, we
%have~$l = 0$.
\smallskip

\noindent
{\em Second case.}
$k = \bbR$.

\noindent
Without~loss of generality, suppose that
$a > b > 0$
and~$a' > b' > 0$.
Then~it will suffice to prove the assertion for representatives of
$e_2$
and~$e_3$,
i.e.,~for
$((x-a)(x-b), u(u-b'))_\bbR$
and
$(x(x-b), (u-a')(u-b'))_\bbR$.
Compare~the proof of Proposition~\ref{dims}.c)~below.

Concerning~$e_2$,
$((x-a)(x-b), u(u-b'))_\bbR = \frac12$
would mean that
$(x-a)(x-b) < 0$
and
$u(u-b') < 0$.
Hence,~$b < x < a$
and
$0 < u < b'$.
But~then
$x(x-a)(x-b) u(u-a')(u-b') < 0$
such that there is no real point
on~$S$
corresponding
to~$(x,u)$.
For~$e_3$,
the argument is~analogous.
}
\eop
\end{prop}

\begin{algo}
\label{testfl}
Let~the parameters
$a,b,a',b' \in \bbZ$,
a Brauer
class~$\alpha \in \Br(S)_2$
as a combination of Hilbert symbols, and a prime
number~$p$
be~given. Then~this algorithm determines the colouring
of~$S(\bbQ_p)$
defined by
$\ev_\alpha\colon S(\bbQ_p) \to \frac12\bbZ/\bbZ$,
for
$S$
the surface given by
$z^2 = x(x-a)(x-b) u(u-a')(u-b')$.

\begin{iii}
%\item
%Check~that
%$\gcd(a,b) = \gcd(a',b') = 1$.
%Otherwise,~output an error message and~terminate.
\item
Calculate~$l := \max (\nu_p(a), \nu_p(b), \nu_p(a-b), \nu_p(a'), \nu_p(b'), \nu_p(a'-b'))$,
the bound established in Lemma~\ref{appr}.
\item
Initialize three lists
$S_0$,
$S_1$,
and
$S_2$,
the first two being empty, the third containing all triples
$(x_0,u_0,p)$
for~$x_0, u_0 \in \{0,\ldots,p-1\}$.
A~triple
$(x_0,u_0,p^e)$
shall represent the subset
$\{(x,u;z) \in S(\bbQ_p) \mid \nu_p(x-x_0) \geq e, \nu_p(u-u_0) \geq e \}$.
\item\label{vier}
Run~through
$S_2$.
For~each element
$(x_0,u_0,p^e)$,
execute, in this order, the following~operations.
\item[$\bullet$ ]
Test~whether the corresponding set is non-empty. Otherwise,~delete~it.
\item[$\bullet$ ]
If~$e \geq l+1$,
$\nu_p(x-\mu) \geq l+1$
for some
$\mu \in \{0,a,b\}$,
and
$\nu_p(u-\nu) \geq l+1$
for a
$\nu \in \{0,a',b'\}$
then move
$(x_0,u_0,p^e)$
to~$S_0$.
\item[$\bullet$ ]
Test~naively, using the elementary properties of the Hilbert symbol, whether the elements in the corresponding set all have the same~evaluation. If~this test succeeds then move
$(x_0,u_0,p^e)$
to~$S_0$
or~$S_1$,~accordingly.
\item[$\bullet$ ]
Otherwise,~replace
$(x_0,u_0,p^e)$
by the
$p^2$
triples
$(x_0 + ip^e, u_0 + jp^e, p^{e+1})$
for
$i, j \in \{0,\ldots,p-1\}$.
\item
If~$S_2$
is empty then output
$S_0$
and
$S_1$
and terminate. Otherwise,~go back to step~\ref{vier}.
\end{iii}
\end{algo}

\begin{ex}
Consider the Kummer
surface~$S$
over~$\bbQ$,
given~by
$$z^2 = x(x-1)(x-25) u(u+25)(u+36) \, .$$
Then~weak approximation is violated
on~$S$.\smallskip

\noindent
{\bf Proof.}
This~is caused by a transcendental Brauer-Manin~obstruction. In~fact,~the matrix~(\ref{szmat})~is
$$M =
\left(\!
\begin{array}{rrrr}
   1 &   25 & 900 &   25 \\
  25 &    1 & -25 & -275 \\
 900 &  -25 &   1 &  -24 \\
  25 & -275 & -24 &    1
\end{array}
\right) 
\widehat{=}
\left(
\begin{array}{rrrr}
   1 &    1 &   1 &    1 \\
   1 &    1 &  -1 &  -11 \\
   1 &   -1 &   1 &   -6 \\
   1 &  -11 &  -6 &    1
\end{array}
\right) ,
$$
having the
kernel~$\langle e_1 \rangle$.
Hence,~there is a transcendental Brauer class
on~$S$,
represented by the quaternion
algebra~$((x-1)(x-25), (u+25)(u+36))$.

Now,~the argument is completely~elementary. For~every
$(x,u;z) \in S(\bbQ_p)$
with
$z \neq 0$,
one has
$$\sum_p ((x-1)(x-25), (u+25)(u+36))_p = 0 \, ,$$
according to the sum formula for the Hilbert~symbol. The~bad primes of the elliptic curves
$y^2 = x(x-1)(x-25)$
and
$y^2 = x(x+25)(x+36)$
are
$2,3,5$,
and~$11$.
Hence,~the sum is actually only over these four~primes.

Our~implementation of Algorithm~\ref{testfl} shows that the local evaluation map is constant at the primes
$2$,
$3$,
and
$11$,
but not
at~$5$.
Hence,~\mbox{$5$-adic}
points such that
$((x-1)(x-25), (u+25)(u+36))_5 = \frac12$
may not be approximated by
\mbox{$\bbQ$-rational}~ones.

Examples~for such
\mbox{$5$-adic}
points are those
with~$(x,u) = (17,5)$.
Indeed,
$17\!\cdot\!(17-1)\!\cdot\!(17-25)\!\cdot\!5\!\cdot\!(5+25)\!\cdot\!(5+36) = -535\,296\!\cdot\!5^2$
is a
\mbox{$5$-adic}~square,
but
$(17-1)\!\cdot\!(17-25) = -128$
is a non-square and
$\nu_5((5+25)\!\cdot\!(5+36)) = 1$
is~odd.
\eop
\end{ex}

\begin{rems}
\begin{iii}
\item
The~constancy of the local evaluation maps at
$3$
and~$11$
and the non-constancy
at~$5$
also follow from the criterion formulated as Theorem~\ref{fund}~below.
\item
In~the colouring obtained
on~$S(\bbQ_5)$,
all the points such that
$x,u \not\equiv 0 \pmod 5$
have colour~zero. This~is rather different from the colourings typically obtained from an algebraic Brauer~class. The~reader should compare the situation described in~\cite{CKS}, where, on the cone over an elliptic curve, three sets of equal sizes~appear.
\end{iii}
\end{rems}

\subsubsection*{Normal Form, Ranks, Asymptotics}

\begin{ttt}[A normal form]
\label{norm}
Let~$k$
be a field,
$a,b,a',b' \in k^*$,
$a \neq b$,
$a' \neq b'$,
and
$S$
be the Kummer surface
$z^2 = x(x-a)(x-b) u(u-a')(u-b')$.
There~are two types of non-trivial Brauer classes
$\alpha \in \Br(S)_2/\Br(k)_2$.\smallskip

\noindent
{\em Type~1.}
$\alpha$~may be expressed by a single Hilbert~symbol.

\noindent
There~are nine cases for the kernel vector
of~$M_{aba'b'}$.
As~seen in Observation \ref{isotw}.i), a suitable translation of
$\smash{\bA^1 \times \bA^1}$
transforms the surface into an isomorphic one with kernel
vector~$e_1$.
Then~$ab, a'b', (-aa') \in k^{*2}$.
Note~that this implies
$(-ba'), (-ab'), (-bb') \in k^{*2}$,~too.
\smallskip

\noindent
{\em Type~2.}
To~express
$\alpha$,
two Hilbert symbols are~necessary.

\noindent
There~are six cases for the kernel vector
of~$M_{aba'b'}$.
A~suitable translation
of~$\smash{\bA^1 \times \bA^1}$
transforms the surface into an isomorphic one with kernel
vector~$e_2+e_3$.
Then~$aa', bb', (a-b)(a'-b') \in k^{*2}$.
\end{ttt}

\begin{coro}
\label{sq}
Let\/~$p$
be a prime number,
$a,b,a',b' \in \bbQ_p^*$,
$a \neq b$,
$a' \neq b'$,
and\/
$S$~be
the Kummer surface\/
$z^2 = x(x-a)(x-b) u(u-a')(u-b')$.
Suppose~that\/
$\nu_p(a) \leq \nu_p(b)$,
$\nu_p(a') \leq \nu_p(b')$,
and\/~$\Br(S)_2/\Br(k)_2 \neq 0$.
Then\/~$\nu_p(aa')$
is~even.\smallskip

\noindent
{\bf Proof.}
{\em
The~assertion is that the expression
$$m := \min(\nu_p(a), \nu_p(b), \nu_p(a-b)) + \min(\nu_p(a'), \nu_p(b'), \nu_p(a'-b'))$$
is even as soon as
$\Br(S)_2/\Br(k)_2 \neq 0$.
As~$m$
is invariant under translations as described in Observation~\ref{isotw}.i), we may suppose that
$e_1 \in \ker M_{aba'b'}$
or
$e_2+e_3 \in \ker M_{aba'b'}$.
In~both cases, the assertion is easily~checked. Note~that either minimum is adopted by at least two of the three~valuations.
\eop
}
\end{coro}

\begin{rems}
\label{nf}
\begin{iii}
\item
Suppose~$k = \bbQ_p$.
Then,~by Observation~\ref{isotw}.ii), we may assume without loss of generality that
$a,b,a',b' \in \bbZ_p$,
$\min(\nu_p(a), \nu_p(b)) = 0$,
and
$\min(\nu_p(a'), \nu_p(b')) = 0,1$.
By~Corollary~\ref{sq}, the assumption that
$M_{aba'b'}$
has a non-trivial kernel ensures that
$\min(\nu_p(a'), \nu_p(b')) = 0$,~too.
\item
Consider the case
that~$k = \bbQ$
and suppose that there is a Brauer class of type~1. According~to i), we may suppose that
$\gcd(a,b) = \gcd(a',b') = 1$.
Hence,~there is a normal form such that
$a > b$,
$a' < b'$,
and~$a,b,(-a'),(-b') \in \bbQ^{*2}$.
Up~to the involution
$(a,b,a',b') \mapsto (-a',-b',-a,-b)$,
this normal form is~unique. Geometrically,~this involution means to interchange the two elliptic curves and to twist both
by~$(-1)$.
\end{iii}
\end{rems}

\begin{prop}
\label{dims}
Let\/
$E \colon y^2 = x(x-a)(x-b)$
and\/
$E' \colon v^2 = u(u-a')(u-b')$
be two elliptic curves over a
field\/~$k$
of characteristic~zero. Suppose~that\/
$a,b,a',b' \in k$
and that\/
$E$
and\/
$E'$
are not isogenous to each~other. Further,~let\/
$S := \Kum(E \times E')$
be the corresponding Kummer~surface.

\begin{abc}
\item
In~all cases,
$\dim \Br(S)_2/\Br(k)_2 \leq 4$
and\/
$\dim \Br(S)_2/\Br(k)_2 \neq 3$.
Further,
$\dim \Br(S)_2/\Br(k)_2 = 4$
is possible only when\/
$(-1)$
is a square
in\/~$k$.
\item
Let\/~$p$
be a prime number
and\/~$k = \bbQ_p$.
If~both\/
$E$
and\/~$E'$
have potential good reduction then\/
$\dim \Br(S)_2/\Br(k)_2$
is~even.
\item
If\/
$k = \bbR$
then\/~$\dim \Br(S)_2/\Br(k)_2 = 2$.
\end{abc}\smallskip

\noindent
{\bf Proof.}
{\em
All~these assertions are consequences of Proposition~\ref{Skoro}. Recall~that
$M_{aba'b'}$
is a matrix with entries in the
$\bbF_{\!2}$-vector
space~$\smash{k^*/k^{*2}}$.\smallskip

\noindent
a)
The~inequality
$\dim \Br(S)_2/\Br(k)_2 \leq 4$
is~clear. Dimension~three would imply that
$M_{aba'b'}$
is of column rank~one. But~this is impossible for a symmetric matrix having zeroes on the~diagonal.
Further,~$\dim \Br(S)_2/\Br(k)_2 = 4$
requires
$M_{aba'b'}$
to be the zero~matrix. In~particular,
$aa'$
and~$(-aa')$
both have to be squares
in~$k$.
This~implies that
$(-1)$
is a square,~too.\smallskip

\noindent
b)
Standard~considerations (cf.~\cite[Proposition~VII.5.5]{Si}) show that the elliptic curve given
by~$y^2 = x(x-\mu)(x-\nu)$
has potential good reduction if and only if
$\mu/\nu \in \bbZ_p^*$
and~$\mu/\nu \not\equiv 1 \pmod p$.
This~implies, in
particular,
that~$p > 2$.

Suppose~that
$\Br(S)_2/\Br(k)_2 \neq 0$
as, otherwise, the assertion is true,~trivially. Then,~by Remark~\ref{nf}.i), we may assume that
$a,b,a-b,a',b',a'-b' \in \bbZ_p^*$.
But~for
\mbox{$p$-adic}
units, being a square
in~$\bbQ_p$
or not is tested by the Legendre~symbol.
$M_{aba'b'}$
is essentially an alternating matrix with entries
in~$\bbF_{\!2}$.
Such~matrices have even~rank.\smallskip

\noindent
c)
Applying~one of the translations
$\smash{\bA^1 \times \bA^1 \to \bA^1 \times \bA^1}$,
$(x, u) \mapsto (x-\mu, u-\nu)$
for
$\mu = 0,a,b$,\,
$\nu=0,a',b'$,
we may assume that
$a > b > 0$
and
$a' > b' > 0$.
Then
$$M_{aba'b'} =
\left(
\begin{array}{cccc}
 + & + & + & - \\
 + & + & + & + \\
 + & + & + & + \\
 - & + & + & +
\end{array}
\right)
$$\vskip-1mm\noindent
has
kernel~$\langle e_2, e_3 \rangle$.
}
\eop
\end{prop}

\begin{rems}[asymptotics]
\begin{iii}
\item
Let~$N>0$.
Then~the number of pairs
$(a,b)$
such that
$a$
and~$b$
are perfect squares,
$a < b$,
and~$a, b-a < N$
is asymptotically
$CN$
for
$C := \frac12 [\log(\sqrt{2}+1) + \sqrt{2} - 1]$.
Indeed,~the Stieltjes integral
$\smash{\int_1^N \!\!\sqrt{x+N} - \sqrt{x} \;d\sqrt{x}}$
has exactly this~behaviour. Assuming~that isogenies are rare, we obtain that the number of surfaces
over~$\bbQ$
with integer coefficients of absolute
value~$\leq\!\! N$
and a
\mbox{$2$-torsion}
Brauer class of type~1 is asymptotically
$\frac12 (6/\pi^2)^2 C^2 N^2 \approx 0.077\,544 N^2$.
\item
On~the other hand, a
\mbox{$2$-torsion}
Brauer class of type~2 yields a
\mbox{$\bbQ$-rational}
point on the intersection of three quadrics
in~$\smash{\bP^6}$.
The~Manin conjecture leads to the naive expectation of a growth of the
type~$cN \log^l \!N$.
\item
The~number of all Kummer surfaces of the form considered and with coefficients up
to~$N$
is~$O(N^4)$.
Thus,~only a very small fraction have a non-trivial
\mbox{$2$-torsion}
Brauer~class.

Even~fewer surfaces should have odd torsion in their Brauer~group. Indeed,~for
\mbox{$l$-torsion,}
one must have
$\smash{\Hom_{\Gal(\overline\bbQ/\bbQ)}(E[l], E'[l]) \neq 0}$~\cite[Propo\-sition~3.3]{SZ}.
Consequently,
$\#E(\bbF_{\!p}) \equiv \#E'(\bbF_{\!p}) \pmod l$
for every prime
$p \neq l$
that is good for both
$E$
and~$E'$.
Based~on this, our computations show that, up to
$N=200$,
no surface has an
\mbox{$l$-torsion}
Brauer class
for~$l \geq 5$.
Further,~at most eight pairs of
\mbox{$j$-invariants}
allow a
\mbox{$3$-torsion}
Brauer~class.
\item
It~is possible
over~$\bbQ$
to have a two dimensional
\mbox{$2$-torsion}
Brauer~group.
For~this, in the normal form of Remark~\ref{nf}.ii), one needs that
$a-b$
and
$b'-a'$
are perfect~squares. Further,~these surfaces have four normal forms instead of two, as there are two Brauer classes of type~1. Corresponding~to a pair of Pythagorean triples, we therefore have two Kummer surfaces, differing from each other by a twist
by~$(-1)$.
The~asymptotics of Pythagorean triples~\cite{BV} shows that there are asymptotically
$\frac4{\pi^4} \log^2 (1+\sqrt{2}) N \approx 0.031\,899 N$
surfaces
over~$\bbQ$
with integer coefficients of absolute
value~$\leq \!N$
and a Brauer group of dimension~two.
\item
Some~actual numbers are listed in the table~below. For~a precise description of the sample, compare paragraph \ref{samp}~below.
\end{iii}

\begin{table}[H]
\begin{center}\vskip-7mm
\caption{Surfaces with a $2$-torsion Brauer class}\vskip2.2mm
{\tiny
\begin{tabular}{|r||r|r|r|}
\hline
bound & dimension~2 &            dimension~1, type~1 & dimension~1, type~2 \\\hline
      &             & (among them algebraic classes) & (none is algebraic) \\\hline\hline
   50 &           0 &             183 (\phantom{0}1) &                 38  \\\hline
  100 &           0 &             766 (\phantom{0}2) &                 98  \\\hline
  200 &           2 &            3049 (\phantom{0}3) &                367  \\\hline
  500 &          12 &           18825 (\phantom{0}4) &               1457  \\\hline
 1000 &          20 &           77249 (\phantom{0}8) &               4398  \\\hline
 2000 &          42 &          305812           (11) &              12052  \\\hline
\end{tabular}
}
\end{center}\vskip-3mm
\end{table}
\end{rems}

\subsubsection*{Trivial evaluation}

\begin{theo}[A criterion]
\label{fund}
Let\/~$p>2$
be a prime number and\/
$0 \neq a,b,a',b' \in \bbZ_p$
such that\/
$a \neq b$
and\/
$a' \neq b'$.
Let\/~$S$
be the Kummer surface, given by\/
$z^2 = x(x-a)(x-b) u(u-a')(u-b')$.
Assume~that\/
$e_1$
is a kernel vector of the
matrix\/~$M_{aba'b'}$
and let\/
$\alpha \in \Br(S)_2$
be the corresponding Brauer~class.

\begin{abc}
\item
Suppose\/~$a \equiv b \not\equiv 0 \pmod p$
or\/~$a' \equiv b' \not\equiv 0 \pmod p$.
Then~the evaluation map\/
$\ev_\alpha\colon S(\bbQ_p) \to \bbQ/\bbZ$
is~constant.
\item
If\/~$a \not\equiv b \pmod p$,
$a' \not\equiv b' \pmod p$,
and not all four numbers are\/
\mbox{$p$-adic}
units then the evaluation map\/
$\ev_\alpha\colon S(\bbQ_p) \to \bbQ/\bbZ$
is~non-constant.
\end{abc}\smallskip

\noindent
{\bf Proof.}
{\em
{\em First step.} Preparations.

\noindent
We~are interested in the Hilbert symbol
$((x-a)(x-b), (u-a')(u-b'))_p$.
Recall~that
$\smash{\frac{a}b, \frac{a'}{b'}, (-bb') \in \bbQ_p^{*2}}$.

A~\mbox{$\bbQ_p$-rational}
point
on~$S$
corresponds to a pair of points on the elliptic curves
$\lambda y^2 = x(x-a)(x-b)$
and
$\lambda v^2 = u(u-a')(u-b')$
for a common value
of~$\lambda$.
The~Hilbert symbol then simplifies
to~$(\lambda x, \lambda u)_p$.\smallskip\pagebreak[3]%

\noindent
{\em Second step.}
$2$-descent.

\noindent
By~$2$-descent, cf.~\cite[Proposition~X.1.4]{Si},
$E\colon Y^2 = X(X-a)(X-b)$
has a point in the square class
of~$x$
if and only if the system
\begin{eqnarray*}
x z_1^2 -  t z_2^2 & = & a \\
x z_1^2 - xt z_3^2 & = & b
\end{eqnarray*}
is~solvable.
Eliminating~$t$,
we obtain
$x^2 z_1^2z_3^2 - x z_1^2z_2^2 = ax z_3^2 - bz_2^2$
or
$$(x z_3^2 - z_2^2)(x z_1^2 - b) = (a-b) xz_3^2 .$$
Division
by~$(-b) xz_3^2$
yields
$\smash{(1 - \frac{z_2^2}{z_3^2} \frac1x)(1 - \frac{z_1^2}b x) = 1 - \frac{a}{b}}$.
In~other words,
$E$
has a point in the square class
of~$x$
if and only if
$\smash{(1 - v^2 x)(1 - w^2 \frac{x}b) = 1 - \frac{a}{b}}$
is~solvable.\smallskip

\noindent
{\em Third step.} Application to the Kummer
surface~$S$.

\noindent
As~$\lambda y^2 = x(x-a)(x-b)$
is equivalent to
$y'^2 = \lambda x (\lambda x - \lambda a)(\lambda x - \lambda b)$
and
$b, (-b')$
are squares, we see that
$S$
has a point with coordinates in the square classes of
$x$
and~$u$
if and only if
\begin{eqnarray*}
\textstyle
(1 - v^2 \lambda x)(1 - w^2 \frac{x}b) & = & \textstyle 1 - \frac{a}{b} \\
\textstyle
(1 - v'^2 \lambda u)(1 - w'^2 \frac{u}{b'}) & = & \textstyle 1 - \frac{a'}{b'}
\end{eqnarray*}
has a solution
$(v,w,v',w',\lambda) \in (\bbQ_p^*)^5$.\smallskip

\noindent
a)
Without~loss of generality, assume that
$a \equiv b \not\equiv 0 \pmod p$
and~$a'/b' \in \bbZ_p$.
Let~$(x,u;z) \in S(\bbQ_p)$
be any point such
that~$z \neq 0$.

Then~Lemma~\ref{solv}.a) shows
$(\smash{\frac{u}{b\hsmash{'}}}, \lambda u)_p = 0$.
Further,~by Lemma~\ref{solv}.c),
$\smash{\frac{x}b}$
or
$\lambda x$
is a square
in~$\bbQ_p$.
In~the case
$\lambda x \in \bbQ_p^{*2}$,
the assertion
$(\lambda x, \lambda u)_p = 0$
is clearly~true.
If~$\smash{\frac{x}b \in \bbQ_p^{*2}}$
then
$\smash{0 = (\frac{u}{b\hsmash{'}}, \lambda u)_p = (-\frac{\lambda}{b\hsmash{'}}, \lambda u)_p = (\frac{\lambda x}{-bb\hsmash{'}}, \lambda u)_p = (\lambda x, \lambda u)_p}$.\smallskip

\noindent
b)
Again~without loss of generality, assume that
$p^2 | a$,
that
$b$
is a unit,
and~$a'/b' \in \bbZ_p$.
We~claim, for
$\lambda = -b$,
there is a point
on~$S$
such that
$x = p$
and~$2 | \nu_p(u)$,
thereby
$\lambda u = (-b)u$
being a non-square.

Indeed,~it is obvious that
$(-bp)(p-a)(p-b) \in \bbQ_p^{*2}$.
Further,~by Hensel's Lemma, it suffices to find a pair
$(U_1, U_2) \in \bbF_{\!p}^* \times \bbF_{\!p}^*$
of non-squares such that
$\smash{(1 - U_1)(1 - U_2) = 1 - \frac{\overline{a}'}{\overline{b}'}}$.
For~this, a counting argument~applies. In~fact, each
$\smash{U_1 \in \bbF_{\!p} \!\setminus\! \{0,1,\frac{\overline{a}'}{\overline{b}'}\}}$
uniquely determines its~partner. As~this set contains
$\smash{\frac{p-1}2}$
non-squares and only
$\smash{\frac{p-5}2}$
squares, the assertion~follows.
}
\eop
\end{theo}

\begin{rems}
\begin{iii}
\item
It~might seem strange to use a descent type argument over a local~field. It~seems to us, however, that a direct argument is neither more elegant nor~shorter.
\item
Using~the descent argument above, we also recover the constancy of the evaluation map in the case of good~reduction. Indeed,~Lemma~\ref{solv}.b) implies that
$\frac{x}b$
or
$\lambda x$
is a square or both have even~valuation. The~first two cases are dealt with as~above.
Otherwise,~$\lambda b$
is a square, hence
$-\frac{\lambda}{b\hsmash{'}}$
is a square, too, and one has to
show~$2 | \nu_p(\lambda u)$.
But~this is implied by Lemma~\ref{solv}.b) when looking at the second~equation.
\end{iii}
\end{rems}

\begin{lem}
\label{solv}
Let\/~$p > 2$
be a prime number
and\/~$A, B \in \bbQ_p^*$,
$Q \in \bbQ_p^{*2}$.
Suppose~that the equation
$(1 - Av^2)(1 - Bw^2) = 1 - Q$
is solvable
in~$\bbQ_p^* \times \bbQ_p^*$.

\begin{abc}
\item
In~all cases,
$(A,B)_p = 0$.
\item
If\/~$Q \in \bbZ_p^*$
then\/
$A \in \bbQ_p^{*2}$,
or\/
$B \in \bbQ_p^{*2}$,
or both,
$A$
and\/~$B$,
are of even valuation.
\item
If\/~$Q \in \bbZ_p^*$
and\/
$Q \equiv 1 \pmod p$
then\/
$A \in \bbQ_p^{*2}$
or\/
$B \in \bbQ_p^{*2}$.
\end{abc}\smallskip

\noindent
{\bf Proof.}
{\em
a)
We~have that
$Av^2 + Bw^2 - AB(vw)^2$
is a~square. When~all three summands are of the same valuation, they must be~units. The~assertion is then clearly~true.
Otherwise,~at most two of the three summands have minimal~valuation. Then~their sum is a square,~too. According~to the definition of the Hilbert symbol~\cite[page~55]{BS},
$(A, B)_p = 0$,
$(A, -AB)_p = 0$,
or~$(B, -AB)_p = 0$.
These~three statements are equivalent to each~other.\smallskip

\noindent
b)
We~have
$\nu_p(1 - Q) \geq 0$.
On~the other hand, if both
$A$
and~$B$
are non-squares then
$\nu_p(1 - Av^2), \nu_p(1 - Bw^2) \leq 0$.
This~implies equality, hence
$Av^2, Bw^2 \in \bbZ_p$.
Both~must be units as
$Av^2 + Bw^2 - AB(vw)^2$
is, by assumption, a square
in~$\bbZ_p^*$.\smallskip

\noindent
c)
If~$A$
and~$B$
were both non-squares then
$\nu_p(1 - A v^2) \leq 0$
and
$\nu_p(1 - B w^2) \leq 0$.
As~$\nu_p(1 - Q) > 0$,
this is a~contradiction.
}
\eop
\end{lem}

\begin{ttt}
Experiments~with Algorithm~\ref{testfl} show surprisingly often that there are non-trivial Brauer classes with trivial
\mbox{$p$-adic}~evaluation.
This~is partially explained by the following~result.
\end{ttt}

\begin{theo}
Let\/~$p>2$
be a prime number and
$a,b,a',b' \in \bbQ_p$\/
be such that\/
$E \colon y^2 = x(x-a)(x-b)$
and\/
$E' \colon v^2 = u(u-a')(u-b')$
are two elliptic~curves. Suppose~that\/
$E$
and\/
$E'$
are not isogenous to each~other. Let,~finally\/
$S$
be the corresponding Kummer~surface. Then

\begin{abc}
\item
if\/~$\dim \Br(S)_2/\Br(\bbQ_p)_2 \geq 2$
then there is a non-zero\/
$\alpha \in \Br(S)_2$
such that\/
$\ev_\alpha$
is the zero~map.
\item
If\/~$\dim \Br(S)_2/\Br(\bbQ_p)_2 = 4$
then the subspace of classes with constant evaluation map is of dimension\/
$4$
when~both,
$E$
and\/~$E'$,
have potential good~reduction. The~dimension
is\/~$3$
when neither curve has potential good reduction
and\/~$2$
in the mixed~case.
\end{abc}\smallskip

\noindent
{\bf Proof.}
{\em
By~Remark~\ref{nf}.i), we may assume without loss of generality that
$a,b \in \bbZ_p$,
not both divisible
by~$p$,
and the same for
$a'$
and~$b'$.
The~case that both curves,
$E$
and~$E'$,
have potential good reduction has thus been treated in Proposition~\ref{CoSk}.\smallskip

\noindent
b)
If~neither curve has potential good reduction then, applying a translation of
$\smash{\bA^1 \times \bA^1}$
as in Observation~\ref{isotw}.i),
we may reduce to the case
$a \equiv b \not\equiv 0 \pmod p$
and
$a' \equiv b' \not\equiv 0 \pmod p$.
Then,~by virtue of Theorem~\ref{fund}.b), the Brauer classes corresponding
to~$\langle e_1, e_2, e_3 \rangle$
have constant evaluation maps, but
$\ev_{e_4}$
is non-constant.

Further,~when only
$E'$
is potentially good, the same arguments show that the Brauer classes corresponding
to~$\langle e_1, e_2 \rangle$
have constant evaluation maps, while those of
$e_3$,
$e_4$,
and~$e_3+e_4$
are~non-constant.\smallskip

\noindent
a)
Only~the case that at least one of the curves
$E$
and
$E'$
does not have potential good reduction requires a~proof. Hence,~we may assume that
$a,b \in \bbZ_p$
and
$a \equiv b \not\equiv 0 \pmod p$.
Then
$ab \in \bbQ_p^{*2}$.

The~upper left
$2 \times 2$-block
of~$M_{aba'b'}$
is~zero.
If~the~block
$\smash{(\atop{\!\!\!\!a'b' \;\;\;\;\;\; aa'}{-aa' \;\; a'(a'-b')})}$
occurring in the lower left has no kernel then the
$2 \times 2$-block 
in the upper right is certainly not the zero~matrix.
Therefore,~$\dim \ker M_{aba'b'} \leq 1$,
a~contradiction. Thus,~there is a Brauer class represented by a vector
from~$\langle e_1, e_2 \rangle$.
By~Theorem~\ref{fund}.a), its evaluation map is~constant.
}
\eop
\end{theo}

\section{A point search algorithm for special Kummer surfaces}

\begin{ttt}
Our surfaces are double covers
of~$\smash{\bP^1 \times \bP^1}$,
given by equations of the form
$$w^2 = f_{ab}(x,y) f_{a'b'}(u,v) \, .$$
Here,~$f_{ab}$
is the binary quartic form
$f_{ab}(x,y) := x y (x - a y) (x - by)$.
Thus,~a point
$\smash{([x:y], [u:v]) \in (\bP^1 \times \bP^1)(\bbQ)}$
leads to a point on the surface if and only if the square classes of
$f_{ab}(x,y)$
and
$f_{a'b'}(u,v)$
coincide or one of them is zero.

We will call the solutions with
$f_{ab}(x,y)$
or
$f_{a'b'}(u,v)$
zero the {\em trivial\/} solutions of the~equation. Obviously,~there is a huge number
of trivial~solutions. Our~aim is to describe an efficient algorithm that searches for non-trivial solutions and does not care about the trivial~ones. In~its simplest version, our algorithm works as~follows.
\end{ttt}

\begin{algo}[Point search]
\label{naiv}
Given two lists
$a_1, \ldots, a_k$
and
$b_1, \ldots, b_k$
and a search
bound~$B$,
this algorithm will simultaneously search for the solutions of all equations of the form
$$w^2 = f_{a_i b_i}(x,y) f_{a_j b_j}(u,v) \, .$$
It will find those with
$|x|,|y|,|u|,|v| \leq B$.

\begin{iii}
\item \label{L_berechnen}
%Compute~a
%bound~$L$
%for the linear factors by~putting
%$$L := B (1 + \max \{|a_i|, |b_i| \mid i = 1,\ldots,k\}) \, .$$
Compute~the
bound~$L := B (1 + \max \{|a_i|, |b_i| \mid i = 1,\ldots,k\})$
for the linear~factors.
\item
Store~the square-free parts of the integers in
$[1, \ldots, L]$
in an
array~$T$.
\item\label{P1_auflisten}
Enumerate~in an iterated loop representatives for all points
$\smash{[x : y] \in \bP^1(\bbQ)}$
with
$x, y \in \bbZ$,
$|x|, |y| \leq B$,
and~$x, y \neq 0$.
\item\label{P1_durchmustern}
For~each point
$[x : y]$
enumerated, execute the operations~below.
\item[$\bullet$ ]
Run a loop over
$i = 1,\ldots,k$
to compute the four linear factors
$x$,
$y$,
$x - a_i y$,
and~$x - b_i y$
of~$f_{a_i,b_i}$.
\item[$\bullet$ ]
Store~the square-free parts of the factors
in~$m_1,\ldots,m_4$.
Use~the table
$T$~here.
\item[$\bullet$ ]
Put
$p_1 := \frac{m_1}{\gcd(m_1,m_2)} \frac{m_2}{\gcd(m_1,m_2)}$,
$p_2 := \frac{m_3}{\gcd(m_3,m_4)} \frac{m_4}{\gcd(m_3,m_4)}$,
and
$$\textstyle p_3 := \frac{p_1}{\gcd(p_1,p_2)} \frac{p_2}{\gcd(p_1,p_2)} \, .$$
Thus,~$p_3$
is a representative of the square class
of~$f_{a_i b_i}(x,y)$.
\item[$\bullet$ ]
Store~the quadruple
$(x,y,i, h(p_3))$
into a~list.
Here,~$h$
is a hash-function.
\item
Sort~the list by the last~component.
\item
Split~the list into~parts. Each~part corresponds to a single value
of~$h(p_3)$.
(At~this point, we have detected all collisions of the hash-function.)
\item
Run~in an iterated loop over all the collisions and check whether
$((x,y,i, h(p_3)), (x',y',i', h(p'_3)))$
corresponds to a solution
$([x : y], [x' : y'])$
of the equation
$w^2 = f_{a_i b_i}(x, y) f_{a_{i'} b_{i'}}(x', y')$.
Output~all the solutions~found.
\end{iii}
\end{algo}

\begin{rems}
\begin{iii}
\item
For practical search bounds
$B$,
the first integer overflow occurs when we multiply
$\smash{\frac{p_1}{\gcd(p_1,p_2)}}$
and
$\smash{\frac{p_2}{\gcd(p_1,p_2)}}$.
But we can think of this reduction
modulo~$2^{64}$
as being a part of our hash-function. Note~that the final check of
$f_{a_i b_i}(x, y) f_{a_{i'} b_{i'}}(x', y')$
being a square can be done without multi-precision integers by inspecting the gcd's of the eight~factors.
\item
One~disadvantage of Algorithm~\ref{naiv} is~obvious. It~requires more memory than is reasonably available by present~standards. We~solved this problem by the introduction of what we call a {\em  multiplicative~paging}. This~is an approach motivated by the simple, additive paging as described in~\cite{EJ1}. In~addition, our memory-optimized point search algorithm is based on the following~observation.
\end{iii}
\end{rems}

\begin{lem}
Let\/~$p$
be a good~prime. Then,~for each pair\/
$(x,y)$
with\/
$\gcd(x,y) = 1$,
at most one of the factors\/
$x$,
$y$,
$(x - ay)$,
and\/
$(x - by)$
is divisible
by\/~$p$.
\eop
\end{lem}

%\noindent
%The~memory-optimized point search algorithm works as~follows.

\begin{algo}[Point search using multivariate paging]
\label{paging}
\begin{iii}
\item
Compute~$L$
and the table of square-class representatives as in Algorithm~\ref{naiv}.
\item
Compute~the upper bound
$C := 2 \max \{ |a_i|, |b_i| \mid i = 1,\ldots,k \}$
for the possibly bad~primes.
\item
Initialize~an array of Booleans of
length~$L$.
Use~the value {\tt false} for the~initialization. We~will call this array the markers of the factors already~treated.
\item\label{pp_schleife}
In~a loop, run over all good primes
below~$L$.
Start~with the biggest prime and stop when the upper
bound~$C$
is reached, i.e., work in {\em decreasing\/}~order. For~each
prime~$p_p$,
execute the steps~below. We~call
$p_p$
the {\em page prime}.
\item[$\bullet$ ]
Run~over all multiples
$m$
of~$p_p$
not exceeding
$L$
and such that the
\mbox{$p_p$-adic}
valuation is~odd. For~each
$m$,
do the~following.
\item[$\bullet\bullet$ ]
Check whether
$m$
is marked as already~treated. In~this case, continue with the
next~$m$.
\item[$\bullet\bullet$ ]
Test~whether
$x$,
$y$,
$x - a_i y$,
or
$x - b_i y$
can represent this~value. Here,~use the constraints
$|x|,|y| \leq B$
and~$i \in \{1,\ldots, k\}$.
\item[$\bullet\bullet$ ]
For each possible representation with
$\gcd(x,y) = 1$,
check whether
$x$,
$y$,
$x - a_i y$,
or
$x - b_i y$
is marked as already~treated. Otherwise,~store the quadruples
$(x,y,i,h(p_3))$
into a~list.
\item[$\bullet\bullet$ ]
Mark~the value of
$m$
as treated and continue with the
next~$m$.
\item[$\bullet$ ]
As~in Algorithm~\ref{naiv}, construct all solutions by inspecting the collisions of the hash-function.
\item
Up~to now, all solutions were found such that
$w$
has at least one prime~factor bigger than the bad-primes bound.
%Note~that
%$p_p$
%is always a prime of good reduction and thus only one of the four factors
%$m_1,m_2,m_3,m_4$
%will be divisible by~it.
To get the remaining ones, use the initial algorithm but skip all
values of
$x, y$
that are marked as treated~factors. Further,~break step \ref{P1_durchmustern} early if
$m_3$
or~$m_4$
is marked as~treated.
\end{iii}
\end{algo}

\begin{rem}
The~last step computes all solutions in smooth~numbers. I.e.,~points such that the square classes of
$f_{ab}(x,y)$
and
$f_{a'b'}(u,v)$
are~smooth. It~is an experimental observation that this step takes only a small fraction of the running time, but gives a large percentage of the~solutions. The~algorithm may easily be modified such that only the solutions in smooth numbers are~found. For~this, the markers for treated factors have to be initialized in an appropriate~way.
\end{rem}

\section{Some experiments}

\begin{ttt}[Colouring by covering---A search for regular colourings]
As~noticed in~\ref{reg}, on various types of surfaces~\cite{Br,EJ3}, the (algebraic) Brauer-Manin~obstruction leads to very regular~colourings. Carrying~this knowledge over to the special Kummer surfaces, given by
$S\colon w^2 = f_4(x,y) g_4 (u,v)$,
one is led to test the~following. For~a
\mbox{$\bbQ$-rational}
point
with~$w \neq 0$,
write
$\lambda w_1^2 = f_4(x,y)$
and~$\lambda w_2^2 = g_4 (u,v)$
and expect the colour to be given by the square class
of~$\lambda$.

For~\mbox{$p$-adic}
points, this defines a colouring with four
($p > 2$),
respectively eight
($p=2$)~colours.
At~the infinite place, the colour is given by the sign
of~$\lambda$.
Motivated~by \cite{Br,EJ3}, we assume that the
\mbox{$p$-adic}
colour of a rational point has a meaning only when
$p$~divides
the conductor of one of the elliptic curves used to
construct~$S$. 
Further,~we restricted ourselves to the square classes of even
\mbox{$p$-adic}
valuation (for the primes of bad~reduction). This~does not exclude all rational points reducing to the singular locus at a bad~prime.

Thus,~we get a colouring of the
\mbox{$\bbQ$-rational}
points with
$2^{k+1}$
colours for a surface with
$k$
relevant odd~primes. Weak~approximation would imply that the colour-map is a~surjection. In~the case of a visible obstruction, we would expect that at most half of the possible colours are in the image of the colour~map.

For~a systematic test, we used the 184 elliptic curves with odd conductor
and~$|a|,|b| < 100$.
This~led to
$16\,836$~surfaces.
The~following table gives an overview of the number of colours that~occurred. %Interestingly,~the test provides no hints for a violation of weak approximation on any of the~surfaces.
 
% Rechenzeit: Suchweite 1000:   1 min 20 sec
%             Suchweite 3000:   11 min 14 sec
%             Suchweite 10000:  109 min
%             Suchweite 10000:  3mon 38 sec (Nur glatte Loesungen)
%             Suchweite 30000:  1111 min = 18.5 Stunden
%             Suchweite 30000:  11 min 42 sec (Nur glatte Loesungen)
%             Suchweite 100000: 275 h
%             Suchweite 100000: 51 min  (Nur glatte Loesungen)

% Weiterhin muss man noch sagen:
% - Bei 6 Primstellen wurde bei
%    7301 von 7409 Flaechen alle 128 Farben getroffen
% - Bei 7 Primstellen wurde bei
%    812 von 1726 Flaechen alle 256 Farben getroffen
% 
% Weiterhin wurden bei 32 von 184 Kurven auf jeder daraus konstruierten 
% Flaeche alle Farben getroffen. Man koennte also jetzt mit den 
% verbleibenden 152 Kurven mit noch groesserer Suchweite weitermachen ... .

\begin{table}[H]
\begin{center}\vskip-7.2mm
{\tiny\tabcolsep3pt
\caption{Regular colourings, numbers of points hit by $\bbQ$-rational points}\vskip0.5mm
\begin{tabular}{|l||r|r|r|r|r|r|r|}
\hline
\#bad primes                 &    2 &       3 &      4 &      5 &        6 &        7 &                 8 \\\hline
\#possible colours           &    8 &      16 &     32 &     64 &      128 &      256 &               512 \\\hline
\#surfaces                   &    4 &     182 &   1678 &   5777 &     7409 &     1726 &                60 \\\hline\hline
\#colours found $h = 1000$   &    8 &  15--16 & 26--32 & 32--64 &  33--127 &  31--157 & 27--\phantom{0}81 \\\hline
\#colours found $h = 3000$   &    8 &      16 & 30--32 & 49--64 &  67--128 &  81--226 &           92--192 \\\hline
\#colours found $h = 10000$  &    8 &      16 &     32 & 57--64 &  93--128 & 142--254 &          207--352 \\
only smooth solutions        &    8 &      16 & 31--32 & 54--64 &  79--128 &  99--236 &          113--197 \\\hline
\#colours found $h = 30000$  &    8 &      16 &     32 & 62--64 & 109--128 & 196--256 &          303--474 \\
only smooth solutions        &    8 &      16 &     32 & 59--64 &  92--128 & 146--253 &          161--300 \\\hline
\#colours found $h = 100000$ &    8 &      16 &     32 &     64 & 121--128 & 232--256 &          387--505 \\
only smooth solutions        &    8 &      16 &     32 & 61--64 & 108--128 & 185--256 &          230--381 \\\hline
\end{tabular}
}
\end{center}\vskip-4mm
\end{table}
\noindent
Our~result is~thus~negative. It~seems that there is no obstruction factoring over such a~colouring. We~expect that one would find
\mbox{$\bbQ$-rational}
points of all colours when further rising the search~bound.

On~one core of an Intel${}^{\text{(R)}}$Core${}^{\text{(TM)}}$2 Duo E8300 processor, the running times were 18.5~hours for search bound
$30\,000$
and 275~hours for search
bound~$100\,000$, but only 51~minutes for smooth solutions with respect to a bad prime bound
of~$200$
and
bound~$100\,000$.
\end{ttt}

\begin{ttt}[Investigating the Brauer-Manin obstruction---A sample]
\label{samp}
We determined all Kummer surfaces of the particular form
$z^2 = x(x-a)(x-b) u(u-a')(u-b')$
that allow coefficients of absolute value 
$\leq \!200$
and have a transcendental
\mbox{$2$-torsion}
Brauer~class.

More~precisely, we determined all
$(a,b,a',b') \in \bbZ^4$
such that
$\gcd(a,b) = 1$,
$\gcd(a',b') = 1$,
$a > b > 0$,
$a - b, b \leq 200$,
as well as
$a' < b' < 0$,
$a' - b', b' \geq -200$
and the matrix
$M_{aba'b'}$
has a non-zero~kernel. We~made sure that
$(a,b,a',b')$
was not listed when
$(-a',-b',-a,-b)$,
$(a,a-b,a',a'-b')$,
or
$(-a',b'-a',-a,b-a)$
was already in the~list. We~ignored the quadruples where
$(a,b)$
and~$(a',b')$
define geometrically isomorphic elliptic~curves.

This~led to 3075 surfaces with a kernel vector of type~1, 367 surfaces with a kernel vector of type~2, and two surfaces with
$\Br(S)_2$
of dimension~two. The~last correspond to
$(25, 9, -169, -25)$
and
$(25, 16, -169, -25)$.
Among~the 3075 surfaces, 26 actually have
$\Br(S)_2 = 0$,
due to a
\mbox{$\bbQ$-isogeny}
between the corresponding elliptic~curves.

The~complete list of these surfaces, the exact equations we worked with, and more details are available on both author's web~pages in a file named {\tt ants\_X\_data.txt}.
\end{ttt}

\begin{ttt}[The BM-relevant primes---$p$-adic point of view]
We~say that a Brauer class
$\alpha \in \Br(S)$
{\em works\/} at a
prime\/~$p$
if the local evaluation map
$\ev_{\alpha,p}$
is non-constant. For~every surface in the sample, using Algorithm~\ref{testfl} and Theorem~\ref{fund}, we determined all the {\em \mbox{BM-rel}\-e\-vant~primes\/}
$p$,
i.e.,~those for which there is a Brauer class working
at~$p$.

For~the two surfaces with
$\Br(S)_2$
of dimension two, the situation is as~follows. In~the case of the coefficient vector
$(25, 9, -169, -25)$,
one Brauer class works at
$2$
and~$13$,
another at
$5$
and~$13$,
and the third at all~three. For~the surface corresponding to
$(25, 16, -169, -25)$,
one Brauer class works at
$3$
and~$13$,
another at
$5$
and~$13$,
and the last at all~three.

Among~the other surfaces, we found no relevant prime six~times, one relevant prime 428~times, two 1577~times, three 1119~times, four 276~times, and five nine~times. Finally,~for
$(196, 75, -361, -169)$,
the Brauer class works at
$2$,
$5$,
$7$,
$11$,
$13$,
and~$19$.

For~three surfaces, it happened that the corresponding elliptic curves were isogenous over a proper extension
of~$\bbQ$.
In~these cases, the Brauer-Manin obstruction is~algebraic. For~two of the surfaces, it worked at one prime and at two for the~last.
\end{ttt}

\begin{ttt}[The BM-relevant primes---$\bbQ$-rational points]
When~the Brauer
class~$\alpha$
works at
$l$~primes
$p_1, \ldots, p_l$,
there are
$2^l$
vectors consisting only of zeroes
and~$\frac12$'s.
By~the Brauer-Manin obstruction, half of them are forbidden as values of
$(\ev_{\alpha,p_1}(x), \ldots, \ev_{\alpha,p_l}(x))$
for
\mbox{$\bbQ$-rational}
points~$x \in S(\bbQ)$.
Using~the point search algorithm~\ref{paging}, we tested whether for every surface in the sample and each vector not forbidden, there is actually a rational~point.

It~turned out that this was indeed the~case. Thus,~no further obstruction becomes visible via this~colouring. However,~in some of the cases, rather high search bounds were~necessary. The~following table shows, for the extreme case of six relevant primes, the number of vectors hit for several search~bounds. Somewhat~surprisingly the smallest solution for each colour was smooth with respect to a bad prime bound
of~$800$.

{\tabcolsep3pt
\begin{table}[H]
\begin{center}\vskip-8mm
\caption{Numbers of vectors in the case $(196, 75, -361, -169)$}\vskip1mm
{\tiny
\begin{tabular}{|r||r|r|r|r|r|r|r|r|r|r|r|}
\hline
bound   & 50 & 100 & 200 & 400 & 800 & 1600 & 3200 & 6400 & 12\,800 & 25\,600 & 50\,000 \\\hline
vectors &  5 &  10 &  14 &  20 &  24 &   26 &   28 &   30 &      31 &      31 &      32 \\\hline
\end{tabular}
}
\end{center}
\end{table}
}\vskip-3.5mm

For~the other surfaces in the sample, lower search bounds were sufficient, but the differences were~enormous. We~summarize our observations in the table~below.

\begin{table}[H]
\begin{center}\vskip-7.7mm
\caption{Search bounds to get all vectors by rational points}\vskip1mm
{\tiny
\begin{tabular}{|r||r|r|r|r|r|r|r|r|r|r|}
\hline
         &            & \multicolumn{9}{c|}{bound $N$ insufficient for}                           \\\hline
\#primes & \#surfaces & $N=50$ & $100$ & $200$ & $400$ & $800$ & $1600$ & $3200$ & $6400$ & $12800$ \\\hline\hline
       2 &       1577 &    190 &    56 &    22 &     - &       &        &        &        &         \\\hline
       3 &       1119 &    555 &   187 &    48 &     1 &     - &        &        &        &         \\\hline
       4 &        262 &    262 &   200 &   127 &    67 &    36 &     24 &     13 &      4 &       - \\\hline
       5 &          9 &      9 &     9 &     8 &     8 &     8 &      5 &      3 &      1 &       - \\\hline
\end{tabular}
}
\end{center}\vskip-2.8mm
\end{table}
\end{ttt}

\begin{rem}
There~is the expectation that the behaviour of the evaluation map
$\ev_{\alpha,p}$
is strongly connected to the type of bad reduction at the
prime~$p$.
For~algebraic Brauer classes, such a connection is well known, cf.\,\cite{EJ3}. In~the transcendental case, there are only partial results, see for example \cite[section~4]{HV}.

For~our examples, the
reductions~$S_p$
are rational surfaces having one or two double~lines.
Further,~$\ev_{\alpha,p}$
is necessarily constant on the set
of~\mbox{$\bbQ$-rational}
points reducing to a smooth~point. The~finer structure seems to be complicated, compare Lemma~\ref{appr}.
\end{rem}


\frenchspacing
\begin{thebibliography}{HVV}
\bibitem[BS]{BS}
Borevich,~Z.\,I.\ and Shafarevich,~I.\,R.: {\em Number theory,} Academic Press, New York-Lon\-don~1966
\bibitem[BV]{BV}
Benito, M.\ and Varona, J.\,L.: {\em Pythagorean triangles with legs less than
$n$,}
J.\ Comput.\ Appl.\ Math.\ {\bf 143}\br(2002)\brr117--126\pagebreak[3]
\bibitem[Br]{Br}
Bright, M.: {\em The Brauer-Manin obstruction on a general diagonal quartic surface,} Acta Arith.\ {\bf 147}\br(2011)\brr291--302
\bibitem[CG]{CG}
Cassels, J.\,W.\,S.\ and Guy, M.\,J.\,T.: {\em On the Hasse principle for cubic surfaces,} Mathematika {\bf 13}\br(1966)\brr111--120
\bibitem[CKS]{CKS}
Colliot-Th\'el\`ene, J.-L., Kanevsky, D., and Sansuc, J.-J.: {\em Arithm\'etique des surfaces cubiques diagonales,} in: Diophantine approximation and transcendence theory (Bonn~1985), Lecture Notes in Math.~1290, Springer, Berlin~1987, 1--108
\bibitem[CS]{CS}
Colliot-Th\'el\`ene, J.-L.\ and Skorobogatov, A.\,N.: {\em Good reduction of the Brauer-Manin obstruction,} arXiv:1006.1972, To appear in: Trans.\ Amer.\ Math.\ Soc.
\bibitem[EJ1]{EJ1}
Elsenhans, A.-S.\ and Jahnel, J.: {\em The Diophantine equation\/
$x^4 + 2 y^4 = z^4 + 4 w^4$,}
Math.\ Comp.\ {\bf 75}\br(2006)\brr935--940
\bibitem[EJ2]{EJ2}
Elsenhans, A.-S.\ and Jahnel, J.: {\em On the Brauer-Manin obstruction for cubic surfaces,} J.\ Combinatorics and Number Theory {\bf 2}\br(2010)\brr107--128
\bibitem[EJ3]{EJ3}
Elsenhans, A.-S.\ and Jahnel, J.: {\em On the quasi group of a cubic surface over a finite field,} J.\ Number Theory {\bf 132}\br(2012)\brr1554--1571
\bibitem[EJ4]{EJ4}
Elsenhans, A.-S.\ and Jahnel, J.: {\em On the order three Brauer classes for cubic surfaces,} Central European Journal of Mathematics {\bf 10}\br(2012)\brr903--926
\bibitem[Ha]{Ha}
Harari, D.: {\em Obstructions de Manin transcendantes,} in: Number theory (Paris 1993--1994), London Math.\ Soc.\ Lecture Note Ser.~235, Cambridge Univ.\ Press, Cambridge~1996, 75--87
\bibitem[HV]{HV}
Hassett, B.\ and V\'arilly-Alvarado, A.: {\em Failure of the Hasse principle on general\/
$K3$~surfaces,}
arXiv:1110.1738
\bibitem[HVV]{HVV}
Hassett, B., V\'arilly-Alvarado, A., and Varilly, P.: {\em Transcendental obstructions to weak approximation on general\/
$K3$~surfaces,} Adv.\ Math.\ {\bf 228}\br(2011)\brr1377--1404
\bibitem[Ie]{Ie}
Ieronymou, E.: {\em Diagonal quartic surfaces and transcendental elements of the Brauer groups,} J.\ Inst.\ Math.\ Jussieu {\bf 9}\br(2010)\brr769--798
\bibitem[ISZ]{ISZ}
Ieronymou, E., Skorobogatov, A., and Zarhin, Yu.: {\em On the Brauer group of diagonal quartic surfaces,} With an appendix by Peter Swinnerton-Dyer, J.\ Lond.\ Math.\ Soc.\ {\bf 83}\br(2011)\brr659--672
\bibitem[Ka]{Ka}
Karzhemanov, I.: {\em One construction of a
$K3$~surface
with the dense set of rational points,} arXiv:1102.1873
\bibitem[LKL]{LKL}
Logan, A., McKinnon, D., and van Luijk, R.: {\em Density of rational points on diagonal quartic surfaces,} Algebra \& Number Theory {\bf 4}\br(2010)\brr 1--20
\bibitem[Ma]{Ma}
Manin, Yu.~I.: {\em Cubic forms,} North-Holland Publishing Co.\ and American Elsevier Publishing Co., Amsterdam-London and New York~1974
\bibitem[Mo]{Mo}
Mordell, L.\,J.: {\em On the conjecture for the rational points on a cubic surface,} J.~London Math.\ Soc.\ {\bf 40}\br(1965)\brr149--158
\bibitem[Pr]{Pr}
Preu, T.: {\em Transcendental Brauer-Manin obstruction for a diagonal quartic surface,} Ph.\,D.\ thesis, Z\"urich~2010
\bibitem[Si]{Si}
Silverman, J.~H.: {\em The arithmetic of elliptic curves,} Graduate Texts in Mathematics~106, Springer, New York~1986
\bibitem[SZ]{SZ}
Skorobogatov, A.\,N.\ and Zarhin, Yu.\,G.: {\em The Brauer group of Kummer surfaces and torsion of elliptic curves,} J.\ reine angew.\ Math.\ {\bf 666}\br(2012)\brr115--140
\bibitem[SSD]{SSD}
Skorobogatov, A.\ and Swinnerton-Dyer, Sir~Peter: {\em
$2$-descent
on elliptic curves and rational points on certain Kummer surfaces,}
Adv.\ Math.\ {\bf 198}\br(2005)\brr448--483
\bibitem[SD]{SD}
Swinnerton-Dyer, Sir~Peter: {\em Two special cubic surfaces,} Mathematika {\bf 9}\br(1962)\brr54--56
\bibitem[Wi]{Wi}
Wittenberg, O.: {\em Transcendental Brauer-Manin obstruction on a pencil of elliptic curves,} in: Arithmetic of higher-dimensional algebraic varieties (Palo Alto~2002), Progr.\ Math.~226, Birkh\"auser, Boston~2004, 259--267
\end{thebibliography}
\end{document}